\documentclass[a4paper,12pt]{amsart}

\usepackage{epsf}

\def\oM{{\overline{\mathcal{M}}}}
\def\C{\mathbb{C}}
\def\Z{\mathbb{Z}}

\def\d{\partial}

\newtheorem{theorem}{Theorem}[section]
\newtheorem{lemma}{Lemma}[section]

\newcommand{\inspic}[1]{\begin{tabular}{c}\epsfbox{#1}\end{tabular}}

\title[A definition of descendants at one point]{A definition of descendants at one point in graph calculus}

\keywords{WDVV equation; genus expansion; dGBV algebras; topological recursion relations.}

\subjclass[2000]{53D45}

\author{Sergei Shadrin}

\address{Department of Mathematics, Stockholm University, Stockholm, SE 106 91, Sweden and Moscow Center for Continuous Mathematical Education, Boljshoi Vlasjevskii Pereulok 11, Moscow, 119002, Russia.}

\email{shadrin@math.su.se and shadrin@mccme.ru}

\thanks{Supported in part by grants RFBR-05-01-01012-a, NSh-1972.2003.1,
NWO-RFBR-047.011.2004.026 (RFBR-05-02-89000-NWO-a),
by the G{\"o}ran Gustafsson foundation, and by
Pierre Deligne's fund based on his 2004 Balzan prize in mathematics.}

\begin{document}

\begin{abstract}
We study the genus expansion of Barannikov-Kon\-tse\-vich solutions of the WDVV equation.
In terms of the related graph calculus we give a definition of descendants at one point
and prove that this definition satisfies the topological recursion relations in genera $0$, $1$, and $2$, string and dilaton equations, and the pull-back formula.
\end{abstract}

\maketitle

\section{Introduction}

\subsection{WDVV equation}
Consider a formal power series $F(T_1,\dots,T_n)$ and a constant non degenerate scalar product $\eta_{ij}$ on the space of variables $T_1,\dots,T_n$. The WDVV equation for $(F,\eta)$ is
\begin{equation*}
\frac{\d ^3 F}{\d T_a \d T_b \d T_i}\eta_{ij}\frac{\d ^3 F}{\d T_j \d T_c \d T_d}
=
\frac{\d ^3 F}{\d T_a \d T_c \d T_i}\eta_{ij}\frac{\d ^3 F}{\d T_j \d T_b \d T_d}.
\end{equation*}
In other words, one can say that $C^{ij}_k=\eta_{lk}\cdot \d ^3 F/\d T_i \d T_j \d T_l$ are the structure constants of an associative commutative algebra.

Solutions of the WDVV equation appear in many natural ways, see~\cite{dz,man1} and references therein. We are interested in two classes of solutions: genus zero Gromov-Witten (GW) invariants~\cite{km} and Barannikov-Kontsevich (BK) construction~\cite{bk}.

The first result that we present in this paper is a very strange property of BK solutions of the WDVV equation. If $(F,\eta)$ is a BK solution of the WDVV equation, then in satisfies the additional equation:
\begin{equation}\label{newE}
\eta_{kl}\frac{\d ^3 F}{\d T_k \d T_l \d T_i}\eta_{ij}\frac{\d ^3 F}{\d T_j \d T_m \d T_n}\eta_{mn}
=
const.
\end{equation}

\subsection{Genus expansion}

Genus zero GW invariants appear naturally to be a part of a bigger formal power series. We can consider GW invariants of higher genera and we can combine them with $\psi$-classes (descendants). Due to the splitting axiom~\cite{km}, any relation among natural strata in the cohomology ring of $\oM_{g,n}$ gives a differential equation for GW invariants.

A genus expansion of BK construction (without descendants) was constructed in~\cite{ls} in terms of tensor expressions associated to graphs. In this paper, we give a partial definition of descendants in terms of graphs. In GW theory our definition corresponds to $\psi$-classes at one point.

It is enough to allow $\psi$-classes at one point in order to fomulate a number of differential equations coming from geometry of the moduli space of curves and therefore universal for GW invariants: topological recursion relations in genera $0$, $1$, and $2$, string and dilaton equations, pull-back formula. We check all these equations for our definition.

There are also Belorousski-Pandharipande relation~\cite{bp} and topological recursion relation in genus $3$~\cite{kl,al2}.
The first one is checked for our definition in~\cite{lss}, and the second one it not yet checked.

\subsection{Further details}

In this paper we work in purely algebraic terms. We start with a vector space $H$ equipped with a structure of so-called
\emph{cyclic Hodge dGBV algebra (cH-algebra)} and we obtain a power series as a contruction of some tensor expressions on this space. Then we describe a technique that allows one to obtain differential equations for this power series.

Thus one of the possible ways to understand this paper is the following. We just offer a formal game with graphs. This game allows one to obtain differential equations, which by a miracle coinside with the equations coming from geometry of the moduli space of curves. In particular, one can use this as a powerful tool allowing to make conjectures on the structure of topological recursion relations in higher genera (cf.~\cite{al1,al2}).

However, all our constructions have strong motivation in geometry, see~\cite[Introduction]{ls}. In particular,
the definition of descendants that we give below comes very natural from the theory of Zwiebach invariants developed in~\cite{ls}. Moreover, using Zwiebach invariants one can define the complete potential that include all descendants.
But study of the complete potential is not yet finished, so we hope to discuss this elsewhere.

\subsection{Acknowledgements}

We are very grateful to A.~Losev for the numerous of helpful remarks and discussions. In this paper we just follow one of his projects. Also we thank B.~Dubrovin, E.~Getzler, S.~Merkulov, and I.~Shneiberg for useful discussions.

\newpage

\section{Construction of the potential}

\subsection{Tensor expressions in terms of graphs}\label{tensors}
We explain a way to encode some tensor expressions over an arbitrary vector space in terms of graphs.

Consider an arbitrary graph (we allow graphs to have leaves). We associate a symmetric $n$-form to each vertex of degree $n$, a symmetric bivector to each egde, and a vector to each leaf. Then we can substitute the tensor product of all vectors in leaves and bivectors in edges into the product of $n$-forms in vertices, distributing the components of tensors in the same way as the corresponding edges and leaves are attached to vertices in the graph. This way we get a number.

Let us study an example:
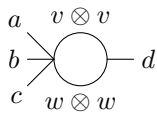
\begin{equation*}
\inspic{arb.1}
\end{equation*}
We assign a $5$-form $x$ to the left vertex of this graph and a $3$-form $y$ to the right vertex. Then the number that we get from this graph is $x(a,b,c,v,w)\cdot y(v,w,d)$.

Note that vectors, bivectors and $n$-forms used in this construction can depend on some variables. Then what we get is not a number, but a function.

\subsection{cH-algebras}\label{cH-algebras}
We recall the definition of cH-alge\-bras~\cite{ls}. A supercommutative associative $\C$-algebra $H$ is called cH-algebra, if there are two odd linear operators
$Q,G_-\colon H\to H$ and an integral $\int\colon H\to\C$ satisfying the following axioms:
\begin{enumerate}
\item $Q^2=G_-^2=QG_-+G_-Q=0$;
\item $H=H_0\oplus H_4$, where $QH_0=G_-H_0=0$ and $H_4$ is represented as a direct sum of subspaces of dimension $4$ generated by $e_\alpha, Qe_\alpha, G_-e_\alpha, QG_-e_\alpha$ for some vectors $e\in H_4$, i.~e. $H_4=\bigoplus_{\alpha}\, \langle e_\alpha, Qe_\alpha, G_-e_\alpha, QG_-e_\alpha \rangle$ (Hodge decomposition);
\item $Q$ is an operator of the first order, it satisfies the Leibniz rule: $Q(ab)=Q(a)b+(-1)^{\tilde a}aQ(b)$
(here and below we denote by $\tilde a$  the parity of $a\in H$);
\item $G_-$ is an operator of the second order, it satisfies the $7$-term relation:
$G_-(abc)=G_-(ab)c+(-1)^{\tilde b(\tilde a+1)}bG_-(ac)
+(-1)^{\tilde a}aG_-(bc)\\
-G_-(a)bc-(-1)^{\tilde a}aG_-(b)c
-(-1)^{\tilde a+\tilde b}abG_-(c).$
\item $G_-$ satisfies the property called $1/12$-axiom: $str(G_-\circ a\cdot)=(1/12)str(G_-(a)\cdot)$ (here $a\cdot$ and $G_-(a)\cdot$ are the operators of multiplication by $a$ and $G_-(a)$ respectively).
\end{enumerate}

We define an operator $G_+\colon H\to H$. We put $G_+H_0=0$; on each subspace $\langle e_\alpha, Qe_\alpha, G_-e_\alpha, QG_-e_\alpha \rangle$, we define $G_+$ as $G_+e_\alpha=G_+G_-e_\alpha=0$, $G_+Qe_\alpha=e_\alpha$, and $G_+QG_-e_\alpha =G_-e_\alpha$. We see that $[G_-,G_+]=0$;
$\Pi_4=[Q,G_+]$ is the projection to $H_4$ along $H_0$; $\Pi_0=\mathrm{Id}-\Pi_4$ is the projection to $H_0$ along $H_4$.

An integral on $H$ is an even linear function $\int\colon H\to\C$. We require
$\int Q(a)b =  (-1)^{\tilde a+1}\int aQ(b)$, $\int G_-(a)b = (-1)^{\tilde a}\int aG_-(b)$, and
$\int G_+(a)b = (-1)^{\tilde a}\int aG_+(b)$. These properties imply that $\int G_-G_+(a)b=\int aG_-G_+(b)$, $\int \Pi_4(a)b=\int a\Pi_4(b)$, and $\int \Pi_0(a)b=\int a\Pi_0(b)$.

We can define a scalar product on $H$: $(a,b)=\int ab.$ We suppose that this scalar product is non-degenerate. Using  scalar product we may turn an operator $A: H \to H$ into bivector that we
denote by $[A]$.

In calculations below, we always use an additional assumption that either the underlying cH-algebra has finite dimension or all tensor expressions that we consider converge.

\subsection{References on cH-algebras}

This structure has appeared first in~\cite{bk} as an axiomatic description of the properties of the holomorphic polyvector fields on Calabi-Yau manifolds. Another example of cH-algebras was found in~\cite{mer1}. In both cases the authors consider only genus zero case, so they didn't include $1/12$-axiom. The full structure of cH-algebra has appeared in~\cite{ls} as an attempt to construct the simplest possible example of Zwiebach invariants.

A detailed discussion of the role of this structure in the theory of Frobenius manifolds one can find in~\cite{man1,man2}. Also the structure of cH-algebra can be obtained as a natural nonlinear generalization of the algebraic structures studied in~\cite{lo,lm,lm2}.

\subsection{Usage of graphs with cH-algebras}\label{usage}

Consider a cH-algebra $H$. There are some standard tensors over $H$, which we associate to elements of graphs below. Here we introduce the notations for these tensors.

We always assign the form
\begin{equation*}
(a_1,\dots,a_n)\mapsto \int a_1\cdot\dots\cdot a_n
\end{equation*}
to a vertex of degree $n$.

There is a collection of bivectors that will arise below at edges: $[G_-G_+]$, $[\Pi_0]$, $[Id]$, $[QG_+]$, $[G_+Q]$, $[G_+]$, and $[G_-]$. In pictures, edges with these bivectors will be denoted by
\begin{equation*}
\inspic{arb.2},\quad \inspic{arb.3},\quad \inspic{arb.4},\quad \inspic{arb.5},\quad \inspic{arb.6},\quad \inspic{arb.7},\quad \inspic{arb.8},
\end{equation*}
respectively. Note that an empty edge corresponding to the bivector $[Id]$ can usually be contracted (if it is not a loop).

The vectors that we will put at leaves depend on some variables. Let $\{e_1,\dots,e_s\}$ be a
homogeneous basis of $H_0$. To each vector $e_i$ we associate formal variables $T_{n,i}$, $n=0,1,2,\dots$, of the same parity as $e_i$. Then we will put at a leaf either the vector $E_0=\sum e_i T_{0,i}$ (denoted by an empty leaf) or the vector $E_n=\sum e_i T_{n,i}$, $n=1,2,\dots$ (then the leaf will be marked by an arrow, the number $n$ will be easily reconstructed from the context).

\subsection{Remark on signs}

Note that $H$ is a $\Z_2$-graded space. In order to get proper signs in our formulas, we always do the following. Suppose we consider a graph of genus $g$. We can choose $g$ edges in such a way that the graph being cut at these edge turns into a tree. To each of these edges we have already assigned a bivector $[A]$ for some operator $A\colon H\to H$. Now we have to put the bivector $[JA]$ instead of the bivector $[A]$, where $J$ is an operator defined by the formula $J\colon a\mapsto (-1)^{\tilde a} a$. Obviously, the result (the number associated to a graph after this procedure) doesn't depend on the choice of these $g$ edges.

For example, consider the following graph:
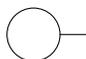
\begin{equation*}
\inspic{arb.9}
\end{equation*}
An empty loop corresponds to the bivector $[Id]$. An empty leaf corresponds to the vector $E_0$. A trivalent vertex corresponds to the $3$-form given by the formula $(a,b,c)\mapsto\int abc$.

If we don't insert $J$, then what we get is just the trace of the operator $a\mapsto E_0\cdot a$. But what we really need is the supertrace of this operator. So, this just will affect on some signs in our calculations.

\subsection{$(g,n)$-Vertices}

In our construction, we consider graphs, whose edges are either marked by $[G_-G_+]$ or by $[Id]$. In the second case, such edges must be loops. All other possible bivectors on edges listed above will appear only in calculations.

Consider a vertex of such graph.
We can split germs of edges attached to this vertex in three groups: $2g$ germs of loops marked by $[Id]$ (the number is even since each loop provide two germs), $k$ germs of edges marked by $[G_-G_+]$ and $m$ germs of leaves.
We say that such vertex is $(g,k+l)$-vertex.

Let us study several examples:
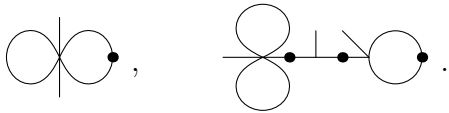
\begin{equation*}
\inspic{arb.10}, \qquad \inspic{arb.11}.
\end{equation*}
The unique vertex of the first graph is $(1,4)$-vertex. The vertices of the second graph are $(2,2)$-, $(0,3)$-, and $(0,4)$-vertex (listed from left to the right).

\subsection{BK solution of WDVV and its genus expansion}

Let $F_0^{sm}=F_0^{sm}(T_{0,1},\dots,T_{0,s})$ be the sum over all trees whose leaves are empty (i.e., marked by $E_0$) and all vertices are $(0,3)$-vertices. We weight each tree by the inverse order of its group of automorphisms.
In other words, we consider trivalent graphs with leaves marked by $E_0$ and edges marked by $[G_-G_+]$:
\begin{equation*}
F_0^{sm}=\frac{1}{6}\inspic{arb.12}+\frac{1}{8}\inspic{arb.13} +\frac{1}{8}\inspic{arb.14} +\dots
\end{equation*}
Let $\eta_{ij}$ be the scalar product of the cH-algebra restricted to $H_0$.

\begin{theorem} \cite{bk,ls} $(F_0^{sm}, \eta)$ is a solution of the WDVV equation.
\end{theorem}

\begin{theorem}\label{strange} $(F_0^{sm}, \eta)$ satisfy Equation~\eqref{newE}.
\end{theorem}

Note that we formulate all equations here and below as if all variables $T_{n,i}$ are even.
It is enough for us since any equation that we prove will first be reformulated in terms of graphs, where all sings are arranged automatically.

Let $F_g^{sm}=F_g^{sm}(T_{0,1},\dots,T_{0,s})$ be the sum over all graphs of genus whose leaves are empty and all vertices are $(0,3)$-vertices. Each graph is weighted by the inverse order of its group of automorphisms. For example,
\begin{align*}
F_1^{sm} &=\frac{1}{2}\inspic{arb.15} +\frac{1}{4}\inspic{arb.16} +\frac{1}{4} \inspic{arb.17} +\dots \\
F_2^{sm} &=\frac{1}{12}\inspic{arb.18}+\frac{1}{8}\inspic{arb.19}
+\frac{1}{8}\inspic{arb.20}+\frac{1}{4}\inspic{arb.21}
+\frac{1}{4}\inspic{arb.22}+\dots
\end{align*}

The formal power series $F^{sm}=\sum_{g\geq 0} F^{sm}_g$ is a natural genus expansion of $F_0^{sm}$ on the small phase space. In GW theory it corresponds to GW invariants in arbitrary genus but without $\psi$-classes.

\begin{theorem} \cite{ls} $(F_1^{sm}, F_0^{sm}, \eta)$ satisfy Getzler elliptic relation~\cite{g1}.
\end{theorem}

\subsection{Descendants at one point}

Now we can describe the formal power series $F$ that corresponds in GW theory to the generating function of GW invariants combined with $\psi$-classes taken only at one point. This formal power series must depend linearly on the variables $T_{n,i}$, where $n\geq 1$ corresponds to the number of $\psi$-classes taken at one point. So, $F$ is the sum over graphs satisfying several conditions:
\begin{enumerate}
\item There is exactly one $(g',m)$-vertex $v_0$ with $3g'-3+m\geq 0$.
\item All other vertices are $(0,3)$-vertices.
\item There is exactly one leaf with arrow.
\item This leaf is marked by $E_n$, $n=3g'-3+m$.
\item This leaf is attached to $v_0$.
\item All other leaves are empty (i.e., marked by $E_0$).
\end{enumerate}
We weight such graph by the inverse order of its group of automorphisms and by $(1/12)^g$.

At the first glance this definition could seem to be a little bit strange. But it is very natural from the point of view of Zwiebach invariants~\cite{ls}. Moreover, as we see below, it is the unique possible definition that satisfies all requirements coming from Gromov-Witten theory.

We can represent $F$ as $\sum_{g\geq 0, n\geq 1} F_{g,n}$ where $F_{g,n}$ is the sum of graphs of total genus $g$ with the special leaf marked by $E_n$. Some examples:
\begin{align*}
F_{0,1}&=\frac{1}{6}\inspic{arb.23}+\frac{1}{4}\inspic{arb.24}+\dots \\
F_{0,2}&=\frac{1}{24}\inspic{arb.25}+\frac{1}{12}\inspic{arb.26}+\dots \\
F_{1,1}&=\frac{1}{24}\inspic{arb.27}+\frac{1}{2}\inspic{arb.28}+\dots \\
F_{1,2}&=\frac{1}{24}\inspic{arb.29}+\frac{1}{48}\inspic{arb.30}+\frac{1}{4}\inspic{arb.31}+\dots \\
F_{2,1}&=\frac{1}{6}\inspic{arb.32}+\frac{1}{4}\inspic{arb.33}+\dots \\
F_{2,2}&=\frac{1}{8}\inspic{arb.34}+\frac{1}{48}\inspic{arb.35}+\dots
\end{align*}

There are some equations for $F^{sm}+F$ that reflect geometry of the moduli space of curves.
By $F_{g,0}$ we denote $F^{sm}_g$, by $F_{g}$ we denote $\sum_{n\geq 0} F_{g,n}$.

\begin{theorem} \textbf{String and dilaton equations.} If $e_1$ is the unit of the underlying cH-algebra, then
\begin{align}
\frac{\d F_{g}}{\d T_{0,1}}&=\sum_{i=1}^s \sum_{n=0}^{\infty} T_{n+1,i}\frac{\d F_{g}}{\d T_{n,i}} + \delta_{g,0}\frac{T_{0,i}\eta_{ij}T_{0,j}}{2}; \label{string}
\\
\frac{\d F_{g,1}}{\d T_{1,1}}&=\sum_{i=1}^s T_{0,i} \frac{\d F_{g,0}}{\d T_{0,i}}+(2g-2)F_{g,0}+\delta_{g,1}\frac{str(\Pi_0)}{24}. \label{dilaton}
\end{align}
\end{theorem}

\begin{theorem} \textbf{Topological recursion relations.} \label{TRR-THM}
\begin{equation}
\frac{\d^3 F_{0,n+1}}{\d T_{n+1,a}\d T_{0,b}\d T_{0,c}}
=\frac{\d^2 F_{0,n}}{\d T_{n,a}\d T_{0,i}}\eta_{ij}\frac{\d^3 F_{0,0}}{\d T_{0,j}\d T_{0,b}\d T_{0,c}};
\end{equation}
\begin{equation}\label{trr1}
\frac{\d F_{1,n+1}}{\d T_{n+1,a}}
=\frac{\d^2 F_{0,n}}{\d T_{n,a}\d T_{0,i}}\eta_{ij}\frac{\d F_{1,0}}{\d T_{0,j}} +\frac{1}{24}
\frac{\d^3 F_{0,n}}{\d T_{n,a}\d T_{0,i}\d T_{0,j}}\eta_{ij};
\end{equation}
\begin{multline}\label{trr2}
\frac{\d F_{2,n+2}}{\d T_{n+2,a}}
=
\frac{\d^2 F_{0,n+1}}{\d T_{n+1,a}\d T_{0,i}}\eta_{ij}\frac{\d F_{2,0}}{\d T_{0,j}}
+
\frac{\d^2 F_{0,n}}{\d T_{n,a}\d T_{0,i}}\eta_{ij}\frac{\d F_{2,1}}{\d T_{1,j}} \\
-
\frac{\d^2 F_{0,n}}{\d T_{n,a}\d T_{0,i}}\eta_{ij}\frac{\d^2 F_{0,0}}{\d T_{0,j}\d T_{0,i'}}\eta_{i'j'}
\frac{\d F_{2,0}}{\d T_{0,j'}} \\
+\frac{7}{10}
\frac{\d^3 F_{0,n}}{\d T_{n,a}\d T_{0,i}\d T_{0,i'}}\eta_{ij}\frac{\d F_{1,0}}{\d T_{0,j}}\eta_{i'j'}
\frac{\d F_{1,0}}{\d T_{0,j'}}\\
+\frac{1}{10}
\frac{\d^3 F_{0,n}}{\d T_{n,a}\d T_{0,i}\d T_{0,i'}}\eta_{ij}\eta_{i'j'}\frac{\d^2 F_{1,0}}{\d T_{0,j}\d T_{0,j'}} \\
-\frac{1}{240}
\frac{\d^2 F_{1,n}}{\d T_{n,a}\d T_{0,i}}\eta_{ij}\frac{\d^3 F_{0,0}}{\d T_{0,j}\d T_{0,i'}\d T_{0,j'}}\eta_{i'j'} \\
+\frac{13}{240}
\frac{\d^4 F_{0,n}}{\d T_{n,a}\d T_{0,i}\d T_{0,j}\d T_{0,i'}}\eta_{ij}\eta_{i'j'}\frac{\d F_{1,0}}{\d T_{0,j'}} \\
+\frac{1}{960}\frac{\d^5 F_{0,n}}{\d T_{n,a}\d T_{0,i}\d T_{0,j}\d T_{0,i'}\d T_{0,j'}}\eta_{ij}\eta_{i'j'}.
\end{multline}
\end{theorem}

We also discuss below the analog of the pull-back formula in terms of graphs.

\subsection{KdV hierarchy}

Consider the one-dimentional cH-algebra generated by the unit $e_1$, with $Q=G_-=0$ and $\int e_1=1$.
Then all graphs contributing to $F^{sm}+F$ have exactly one vertex. It is easy to see that in this case
\begin{equation*}
F^{sm}+F= \frac{T_{0,1}^3}{6} + \sum_{n=1}^{\infty} \frac{T_{n,1}T_{0,1}^{n+2}}{(n+2)!}+
\sum_{g=1}^{\infty}\sum_{n=0}^{\infty}\frac{T_{3g+n-2,1}T_{0,1}^{n}}{g!\cdot 24^g\cdot n!}.
\end{equation*}
This is exactly the part of the string solution of the KdV hierarchy corresponding to the intersection numbers of $\psi$-classes at one point.

\subsection{Organization of the paper}

The rest of the paper is organized like as follows. First, we prove string and dilaton equations. They are almost obvious and don't require any complicated technique. Then we discuss some technical lemmas closed to Tillmann's theorem in~\cite{t}, and we apply this new technique to the analog of the pull-back formula in terms of graphs. Then, using the pull-back formula and an argument shared in~\cite{ls}, we prove topological recursion relations.

And finally we prove the new relation for $F_{0,0}$, which has no relation to the definition of descendants.


\section{String and dilaton equations}

\begin{proof}[Proof of the string equation]
In terms of graph calculus, $\d F_{g,n}/\d T_{0,1}$ is the sum over the same graphs as $F_{g,n}$, but one of simple (empty) leaves is marked not by $E_0$ but by $e_1$. We must also recalculate the coefficiets of graphs taking into account that now the this leaf is preserved by automorphisms of graphs.

In the simplest case we obtain
\begin{equation*}
\frac{1}{2}\inspic{arb.36}.
\end{equation*}
Since $e_1$ is the unit of algebra, this expression is equal to $T_{0,i}\eta_{ij}T_{0,j}/2$. It is exactly the last term in the string equation~\eqref{string}

Now consider a generic graph. Let the leaf marked by $e_1$ be attached to a $(0,3)$-vertex. There are two possible local pictures (in addition to the case studied above):
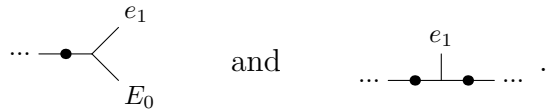
\begin{equation*}
\inspic{arb.37} \qquad \mathrm{and} \qquad \inspic{arb.38}.
\end{equation*}
In the first picture, this piece of graph can be substituted by the vector $G_-G_+(E_0\cdot e_1)=G_-G_+(E_0)=0$. In the second picture, this piece of graph can be substituted by the bivector $[G_-G_+\circ e_1\cdot\circ G_-G_+]=[G_-G_+G_-G_+]=0$ ($e_1\cdot$ is the operator of multiplication by $e_1$). Therefore, in both cases such graphs contribute zero.

So, the leaf marked by $e_1$ can be attached only to the same vertex as the leaf marked by arrow. In particular, if $n=0$,
the unique nonzero contribution to $\d F_{g,n}/\d T_{0,1}$ is the simplest case studied above. Then, since $e_1$ is the unit of the algebra, we can erase it. Thus we obtain a graph, contributing to $F_{g,n-1}$. Moreover, one can easily see that we obtain this graph with the proper combinatorial coefficient.
\end{proof}

\begin{proof}[Proof of the dilaton equation]
First, we note that
\begin{equation*}
\frac{\d}{\d T_{1,1}}\inspic{arb.27}=str(\Pi_0).
\end{equation*}
This gives the exceptional summand in the dilaton equation~\eqref{dilaton}

In all other cases $\d F_{g,1}/\d T_{1,1}$ is represented as a sum over the same graphs as $F_{g,1}$, but the leaf with arrow is marked by $e_1$. Since $e_1$ is the unit of the algebra, we can erase the leaf with arrow, and we obtain a graph contributing to $F_{g,0}$. Since the arrow could be attached to an arbitrary vertex of a graph contributing to $F_{g,0}$,
we get each such graph with the additional coefficient equal to the number of its vertices. This is equal to $2g-2+l$, where $l$ is the number of leaves. The right hand side of dilaton equation is exactly the operator of multiplication of each graph by this coefficient.
\end{proof}


\section{Two technical lemmas}

In this section, we present some technical lemmas. The first one is very closed to the theorem of U.~Tillmann~\cite{t}.
Note, that graphs in this section have a little bit deifferent meaning, then in all other parts of this paper. All graph here will appear later in calculations as parts of bigger graphs. So, empty leaves in graphs of the first lemma don't marked by any vectors, but a graph with $1$ ($2$) empty leaves is considered as a special vector (bivector) itself.

\begin{lemma}\label{tlemma}
The vectors and bivectors defined by the pictures below are equal to zero:
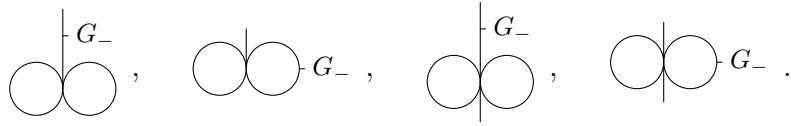
\begin{equation*}
\inspic{arb.39}, \quad \inspic{arb.40}, \quad \inspic{arb.41}, \quad \inspic{arb.42}.
\end{equation*}
\end{lemma}

\begin{lemma}\label{alemma}
For any vectors $A_0,A_1,\dots,A_k$, $k\geq 2$, we have:
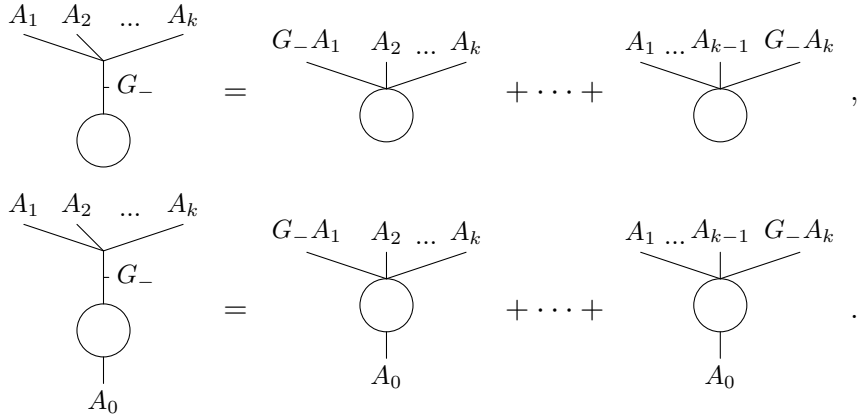
\begin{align*}
\inspic{arb.43} & = \inspic{arb.44}+\dots+\inspic{arb.45}, \\
\inspic{arb.46} & = \inspic{arb.47}+\dots+\inspic{arb.48}.
\end{align*}
\end{lemma}

\begin{proof}[Proof of both lemmas] Both lemmas easily follows from the properties of $G_-$. For instance,
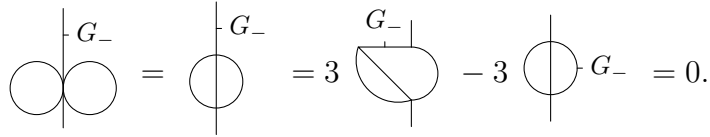
\begin{equation*}
\inspic{arb.41}=\inspic{arb.49}=3\inspic{arb.50}-3\inspic{arb.51}=0.
\end{equation*}
The first equality is a redrawing, the second equality is the special case of th $7$-term relation, and the third equality is again a redrawing.

Another example. Consider the special case of the $7$-term relation:
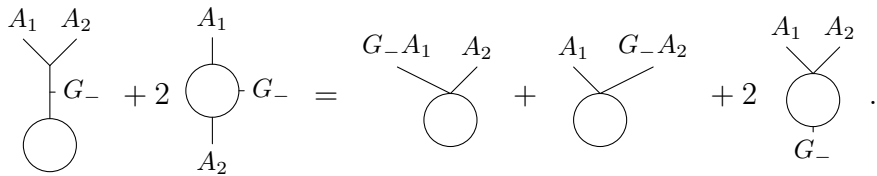
\begin{equation*}
\inspic{arb.52} +2\inspic{arb.53} = \inspic{arb.54}+\inspic{arb.55}+2\inspic{arb.56}.
\end{equation*}
The last summands in both sides of this equation coinside. Thus we obtain the first statement of Lemma~\ref{alemma} for $k=2$.

The other $3$ statements of Lemma~\ref{tlemma} and statements of Lemma~\ref{alemma} for an arbitrary $k$ are proved by the very similar one-string calculations.
\end{proof}


\section{Pull-back of descendants}

\subsection{Pull-back formula on the moduli space of curves}
We consider the moduli space of curves $\oM_{g,r+1}$. We denote by $\pi\colon\oM_{g,r+1}\to\oM_{g,r}$ the projection forgetting the last marked point. Then there is a formula relating $\psi$-classes on $\oM_{g,r+1}$ and pull-backs of $\psi$-classes on $\oM_{g,r}$. We have:
\begin{equation}\label{pbf}
\psi_1^n=\pi^*\psi_1^n+D\cdot \psi_1^{n-1}.
\end{equation}
Here we denote by $D$ the cohomology class of the divisor in $\oM_{g,r+1}$, whose generic point is represented by a two-component curve such that one component has genus $0$ and contains the first and the last marked points and the other component has genus $g$ and contains all other marked points.

\subsection{Pull-back formula in terms of graphs}\label{description}
We give an interpretation of Equation~\eqref{pbf} in terms of graphs. The left hand side of Equation~(\ref{pbf}) is obvious: we consider the sum of graphs with $r$ empty leaves contributing to $F_{g,n}$. Since we don't distinguish the $(r+1)$th leaf on the left hand side, we must multiply the combinatorial coefficient of each graph by $r$.

Now we describe the right hand side of Equation~\eqref{pbf}. The first summand is obtained by the procedure studied in~\cite{ls}. We take a graph with $r-1$ empty leaves contributing to $F_{g,n}$. Then we make one change in this graph in all possible ways. We change either an edge marked by thick black point or a leaf using the rules:
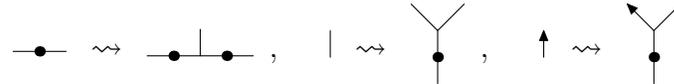
\begin{equation*}
\inspic{arb.57}\leadsto \inspic{arb.58},\quad
\inspic{arb.59}\leadsto \inspic{arb.60},\quad
\inspic{arb.61}\leadsto \inspic{arb.62}
\end{equation*}

In order to obtain the second summand, we take a graph with $r-1$ empty leaves contributing to $F_{g,n-1}$ and substitute
$E_{n-1}$ with $\Pi_0(E_nE_0)$. In the case $n=1$, we do only one change but in all possible ways. Graphically, the rules are:
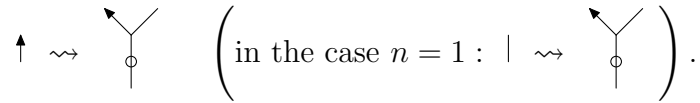
\begin{equation*}
\inspic{arb.63}\leadsto \inspic{arb.64}\quad\left(\mathrm{in\ the\ case\ }n=1:
\inspic{arb.65}\leadsto \inspic{arb.66} \right).
\end{equation*}

\subsection{Example in genus $1$}
We consider projection $\pi\colon\oM_{1,2}\to\oM_{1,1}$. The pull-back formula for the class $\psi_1$ is $\psi_1=\pi^*\psi_1+D$. Its interpretation in terms of graphs:
\begin{equation*}
\frac{1}{2}\inspic{arb.67}=\frac{1}{24}\inspic{arb.68}+\frac{1}{2}\inspic{arb.69}.
\end{equation*}
Let us prove this formula. It follows from
\begin{equation*}
\inspic{arb.69}=\inspic{arb.67}-\inspic{arb.70}=\inspic{arb.67}-\frac{1}{12}\inspic{arb.68}.
\end{equation*}
The first equality is a corollary of $\Pi_0=Id-QG_+-G_+Q$ and the Leibniz rule for $Q$ (see~\cite{ls,lss}, where this standard step is dicussed in all details), the second equality is a corollary of $1/12$-axiom.

\subsection{Example in genus $2$} We consider projection $\pi\colon\oM_{1,2}\to\oM_{1,1}$. The interpretation of the formula $\psi_1^4=\pi^*\psi_1^4+D\cdot \psi_1^{3}$ in terms of graphs is
\begin{equation*}
\frac{1}{48}\inspic{arb.71} = \frac{1}{8\cdot 12^2} \inspic{arb.72} + \frac{1}{48}  \inspic{arb.73}.
\end{equation*}

Let us prove this formula. From Lemma~\ref{tlemma}, it follows that the first summand in the right hand side is equal to zero. Then the pull-back folrmula follows from
\begin{equation*}
\inspic{arb.73}=\inspic{arb.71}-\inspic{arb.74}=\inspic{arb.71}.
\end{equation*}
Here the first equality is a corollary of $\Pi_0=Id-QG_+-G_+Q$ and the Leibniz rule for $Q$, the second equality is a corollary of Lemma~\ref{tlemma}.

\subsection{Proof of the pull-back formula}

\begin{theorem}\label{pbfT}
Descendants in our construction satisfy the pull-back formula.
\end{theorem}

\begin{proof}
First of all, return to the description of the pull-back formula given in Section~\ref{description}.
The summand corresponding to $\pi^*\psi_1^n$ splits into sum $A'+A''$, where $A'$ includes graphs obtained by change at an edge or at an empty leaf, and $A''$ includes graphs obtained by change at the leaf with arrow.

Consider the summand corresponding to $D\cdot \psi_1^{n-1}$. Since $\Pi_0=Id-G_+Q-QG_+$, we can represent each graph in this summand as a sum of three graphs. Then we collect the graphs with $Id$, $-QG_+$, and $-G_+Q$, and denote their sums by $B'$, $B''$, and $B'''$.

We redraw the graphs in $B'$ contructing the edge marked by $[Id]$. Obviously, the sum $A'+B'$ consists of exactly the same graphs as the left hand side of the pull-back formula, and one can easily check that the combinatorial coefficients at the same graphs coinside. So, we must prove that $A''+B''+B'''=0$.

Note that $-G_+Q(E_nE_0)=0$ (apply the Leibniz rule and note that $QE_i=0$, $i\geq 0$). Therefore, $B'''=0$

Consider graphs in $A''$. If such graph has a $(g,n)$-vertex with $g\geq 2$, then it is equal to zero (Lemma~\ref{tlemma}).
Consider graphs in $B''$. Using the Leibniz rule for $Q$, we turn each graph into the sum of graphs of the same type, but with $-[G_+]$ instead of $-[QG_+]$, and with $[G_-]$ instead of one of $[G_-G_+]$. There are several possible cases:
\begin{enumerate}
\item Edges with $-[G_+]$ and $[G_-]$ are attached to different vertices. Using the $7$-term relation and $1/12$-axiom one can show that the sum of all these graphs is equal to zero.
\item Edges with $-[G_+]$ and $[G_-]$ are attached to the same $(g,n)$-vertex, $g\geq 2$. Each such graph is equal to zero (Lemma~\ref{tlemma}).
\item Edges with $-[G_+]$ and $[G_-]$ are attached to the same $(1,n)$-vertex. Each such graph is equal to zero (apply Lemma~\ref{alemma} two times).
\item Edges with $-[G_+]$ and $[G_-]$ are attached to the same $(0,n)$-vertex and $[G_-]$ is on the loop. Applying the $1/12$-axiom and Lemma~\ref{alemma}, we obtain exactly the graphs in $A''$ that have a $(1,n-2)$-vertex.
\item Edges with $-[G_+]$ and $[G_-]$ are attached to the same $(0,n)$-vertex and $[G_-]$ is not on the loop. Applying the $7$-term relation, we obtain exactly the graphs in $A''$ that have a $(0,n+1)$-vertex.
\end{enumerate}
One can easily check that the coefficients of graphs in $A''$ and of the same graphs obtained in $B''$ are opposite. So, $A''+B''=0$.
\end{proof}


\section{Topological recursion relations}

\begin{proof}[Proof of Theorem~\ref{TRR-THM}]
In GW theory, the topological recursion relations appear in the following way. Suppose we have a relation in cohomology of $\oM_{g,k}$ expressing $\psi_1^m$ in terms of boundary strata. Using the splitting axiom, we can understand this relation as a differential equation for the GW potential proven when all parameters are set to zero. Then we can consider the pull-backs of this relation under the projections $\pi\colon\oM_{g,k+n}\to\oM_{g,k}$. This will give
a differential equation for the GW potential proven for arbitrary parameters. But this differential equation differs from the initial one; we must add several terms related to the pull-back formula~\eqref{pbf}. This way we get the first summand in~\eqref{trr1} and the first three summands in~\eqref{trr2}. However, the point is that in GW theory we are to check a TRR-type relation only once, in the cohomology of $\oM_{g,k}$ with the smalest possible $k$.

The advantage of our graph calculus is that in order to prove a topological recursion relation, we can follow the same way as in GW theory. Indeed, it was shown in~\cite{ls} that $\pi^*$ keeps equalities in graphs. In this paper, we proved the pull-back formula in terms of graphs (Theorem~\ref{pbfT}). So, it is sufficient to prove the topological recursion relations when all parameters are set to zero. The difference between GW theory and graph calculus is the following. In GW theory a TRR for $\psi_1^k$ implies the same TRR for $\psi_1^{k+1}$. In graph calculus, there is no analogue of multiplication in cohomology, so we must check the simplest case of each TRR for \emph{all} powers of $\psi$-class.

So, calculations below complete the proof of Theorem~\ref{TRR-THM}.
\end{proof}

\subsection{Simple cases}

In genus $0$ and genus $1$ we just check that
\begin{equation*}
\inspic{arb.75}=\inspic{arb.76} \quad \mathrm{and} \quad
\frac{1}{24}\inspic{arb.77}=\frac{1}{12}\cdot\frac{1}{2}\inspic{arb.78}.
\end{equation*}
(for genus 1, the coefficients $1/24$ and $1/2$ are the weights determined by our rules, the coefficient $1/12$ comes from Equation~\eqref{trr1}). Obviously, this is true.

In genus $2$, if the power of $\psi$-class is $\geq 4$, we just check that
\begin{equation*}
\frac{1}{8\cdot 12^2}\inspic{arb.79}=
-\frac{1}{120}\cdot\frac{1}{48}\inspic{arb.81}
+\frac{1}{120}\cdot\frac{1}{8}\inspic{arb.80}.
\end{equation*}
(the coefficients $1/(8\cdot 12^2)$, $1/48$, $1/8$ are the weights determined by our rules, the coefficients $-1/120$ and $1/120$ come from Equation~(\ref{trr2})). Since all pictures are equal, we are just to check that
$1/8\cdot 12^2=-1/120\cdot 1/48 + 1/120\cdot 1/8$. Obviously, this is true.

\subsection{Genus $2$ expression for $\psi_1^2$}

For convenience, we rewrite the corresponding equation as an expression for $\psi_1^2$ in boundary strata of $\oM_{2,1}$:
\begin{equation}\label{TRR2sim}
\psi_1^2=\frac{7}{5}\inspic{arb.82}+\frac{1}{5}\inspic{arb.83}
-\frac{1}{120}\inspic{arb.84}+\frac{13}{120}\inspic{arb.85}+\frac{1}{120}\inspic{arb.86}.
\end{equation}
In these pictures, a vertex marked by $1$ corresponds to a genus $1$ curve, a simple vertex corresponds to a genus $0$ curve, an edge corresponds to a point of intersection, and a leaf corresponds to the marked point. We note that these pictures have completely different meaning then all other pictures in this paper.

To each stratum we associate a differential monomial in $F^{sm}+F$, and below we list the constant term of this monomials in terms of our usual graphs:
\begin{align*}
\inspic{arb.82} & = \frac{1}{8}\inspic{arb.118} \\
\inspic{arb.83} & = \frac{1}{4}\inspic{arb.119} + \frac{1}{4}\inspic{arb.120} \\
\inspic{arb.84} & = \frac{1}{4}\inspic{arb.121} + \frac{1}{4}\inspic{arb.122} \\
\inspic{arb.85} & = \frac{1}{2}\inspic{arb.123} + \frac{1}{4}\inspic{arb.124} \\
\inspic{arb.86} & = \frac{1}{2}\inspic{arb.125} + \frac{1}{8}\inspic{arb.126} \\
& \phantom{=\ } + \inspic{arb.127} + \frac{1}{4}\inspic{arb.128}
\end{align*}
Using the standard argument in~\cite{ls,lss}, we express all these graphs in terms of the following $5$ pictures:
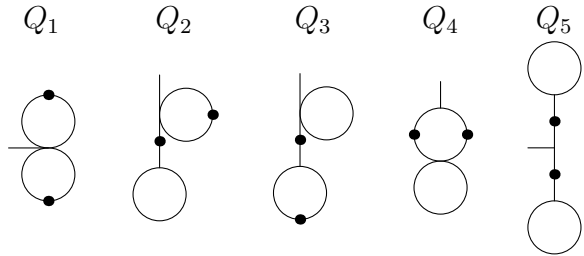
\begin{equation*}
\begin{array}{ccccc}
Q_1 & Q_2 & Q_3 & Q_4 & Q_5 \\
\inspic{arb.113} & \inspic{arb.114} & \inspic{arb.115} & \inspic{arb.116} &
\inspic{arb.117}
\end{array}
\end{equation*}
It is rather hard but straightforward calculation. The result is:
\begin{align*}
\inspic{arb.82} & = \frac{1}{8} Q_1 - \frac{1}{48} Q_2 + \frac{1}{1152} Q_5 \\
\inspic{arb.83} & = -\frac{1}{4} Q_1 + \frac{1}{4} Q_3 + \frac{1}{48} Q_4 \\
\inspic{arb.84} & = -\frac{1}{4} Q_2 + \frac{1}{4} Q_3 + \frac{1}{4} Q_4 \\
\inspic{arb.85} & = \frac{1}{4} Q_2 - \frac{1}{4} Q_3 - \frac{1}{48} Q_5 \\
\inspic{arb.86} & = -\frac{1}{4} Q_4 + \frac{1}{8} Q_5
\end{align*}
If we substitute these expressions in Equation~\eqref{TRR2sim}, then we obtain $\psi_1^2=(1/8)Q_1+(1/48)Q_3$. We see that it is exactly our definition of $\psi_1^2$!

\subsection{Genus $2$ expression for $\psi_1^3$}
For convenience, we rewrite the corresponding equation as an expression for $\psi_1^2$ in boundary strata of $\oM_{2,1}$:
\begin{equation}\label{TRRpsi3}
\psi_1^3=-\frac{1}{120}\psi_1\cdot\inspic{arb.84}+\frac{13}{120}\psi_1\cdot\inspic{arb.85}
+\frac{1}{120}\psi_1\cdot\inspic{arb.86}.
\end{equation}
We list the constant terms of the differential monomials of these strata:
\begin{align*}
\psi_1\cdot\inspic{arb.84} &= \frac{1}{4}\inspic{arb.131},\\
\psi_1\cdot\inspic{arb.85} &= \frac{1}{4}\inspic{arb.132},\\
\psi_1\cdot\inspic{arb.86} &= \frac{1}{4}\inspic{arb.133}+ \frac{1}{2}\inspic{arb.134}.
\end{align*}
Using the standard argument in~\cite{ls,lss}, we express all these graphs in terms of the following $2$ pictures:
\begin{equation*}
P_1=\inspic{arb.129}\qquad \mathrm{and}\qquad P_2=\inspic{arb.130}.
\end{equation*}
We have:
\begin{align*}
& \psi_1\cdot\inspic{arb.84} =\psi_1\cdot\inspic{arb.85} = \frac{1}{4}P_1-\frac{1}{48}P_2, \\
& \psi_1\cdot\inspic{arb.86} = -\frac{1}{2}P_1+\frac{1}{4}P_2.
\end{align*}
If we substitute these expressions in Equation~\eqref{TRRpsi3}, then we obtain $\psi_1^3=(1/48)P_1$. It is exactly our definition of $\psi_1^3$.


\section{New relation}

\begin{proof}[Proof of Theorem~\ref{strange}]
First, let us rewrite Equation~\eqref{newE} as
\begin{equation*}\label{diffnew}
\eta_{i'j'}
\frac{\d^4 F_{0,0}}{\d T_{0,a} \d T_{0,i'} \d T_{0,j'} \d T_{0,i}}
\eta_{ij}
\frac{\d^3 F_{0,0}}{\d T_{0,j} \d T_{0,i''} \d T_{0,j''}}
\eta_{i''j''}
=0
\end{equation*}
We see that this equation is of the same type as equations studied above, so it is enough to prove it  when all parameters are set to zero. In terms of our graphs, we must prove that
\begin{equation}\label{picnew}
\frac{1}{4}\inspic{arb.135}+\frac{1}{2}\inspic{arb.136}=0
\end{equation}
Using (twice) $\Pi_0=Id-QG_+-G_+Q$ and the Leibniz rule for $Q$, one can show that
\begin{align*}
\inspic{arb.135} &=\inspic{arb.137}-\inspic{arb.138},\\
\inspic{arb.136} &=\inspic{arb.139}-\inspic{arb.140}.
\end{align*}
The $7$-term relation for $G_-$ implies that
\begin{equation*}
\inspic{arb.138}+2\inspic{arb.140}=2\inspic{arb.139}.
\end{equation*}
Therefore, the left hand side of Equation~\eqref{picnew} is equal to
\begin{equation*}
\frac{1}{4}\inspic{arb.137}.
\end{equation*}
Using $\Pi_0=Id-QG_+-G_+Q$ and the Leibniz rule for $Q$, we have:
\begin{equation*}
\inspic{arb.137}=\inspic{arb.141}-\inspic{arb.142}.
\end{equation*}
From Lemma~\ref{alemma} is follows that this difference is equal to zero.
\end{proof}



\begin{thebibliography}{00}

\bibitem{al1} D. Arcara, Y.-P. Lee, Tautological equation in $\oM_{3,1}$ via invariance conjectures, math.AG/0503184.

\bibitem{al2} D. Arcara, Y.-P. Lee, Tautological equations in genus $2$ via invariance conjectures, math.AG/0502488.

\bibitem{bk} S.~Barannikov, M.~Kontsevich, Frobenius manifolds and formality of Lie
algebras of polyvector fields, Internat. Math. Res. Notices \textbf{1998}, no. 4, 201--215.

\bibitem{bp} P.~Belorousski, R.~Pandharipande, A descendent relation in genus 2,
Ann. Scuola Norm. Sup. Pisa Cl. Sci. (4) \textbf{29} (2000), no. 1, 171--191.

\bibitem{dz} B.~Dubrovin, Y.~Zhang, Normal forms of hierarchies of integrable PDEs, Frobenius manifolds and Gromov-Witten invariants, arXiv: math.DG/0108160.

\bibitem{g3} E.~Getzler, Topological recursion relations in genus $2$, Integrable Systems and Algebraic Geometry (Kobe/Kyoto, 1997), World Scientific, River Edge, NJ, 1998, pp. 73--106.

\bibitem{g1} E.~Getzler, Intersection theory on $\oM_{1,4}$ and elliptic Gromov-Witten invariants, J. Amer. Math. Soc. \textbf{10} (1997), no. 4, 973--998.

\bibitem{kl} T. Kimura, X. Liu, A genus-3 topological recursion relation, math/0502457.

\bibitem{k} M.~Kontsevich, Intersection theory on the moduli space of curves and the matrix Airy function, Comm. Math. Phys. \textbf{147} (1992), 1--23.

\bibitem{km} M.~Kontsevich, Yu.~Manin, Gromov-Witten classes, quantum cohomology, and enumerative geometry,
Comm. Math. Phys. \textbf{164} (1994), no. 3, 525--562.

\bibitem{lo} A.~Losev, Hodge strings and elements of K. Saito's theory of primitive form,
Topological field theory, primitive forms and related topics (Kyoto, 1996), 305--335, Progr. Math., 160,
Birkhäuser Boston, Boston, MA, 1998.

\bibitem{lm} A.~Losev, Yu.~Manin, Extended modular operad, Frobenius manifolds, 181--211, Aspects Math., E36, Vieweg, Wiesbaden, 2004.

\bibitem{lm2} A.~Losev, Y.~Manin, New moduli spaces of pointed curves and pencils of flat connections, Michigan Math. J.  \textbf{48} (2000), 443--472.

\bibitem{ls} A.~S.~Losev, S.~V.~Shadrin, From Zwiebach invariants to Getzler relation, math.QA/0506039.

\bibitem{man1} Yu.~I.~Manin, Frobenius Manifolds, Quantum Cohomology, and Moduli Spaces, Providence, Rhode Island, 2000.

\bibitem{man2} Yu.~I.~Manin, Three constructions of Frobenius manifolds: a comparative study, Asian J. Math. \textbf{3} (1999), no. 1, 179--220.

\bibitem{mer1} S.~A.~Merkulov, Formality of canonical symplectic complexes and Frobenius
manifolds, Internat. Math. Res. Notices 1998, no. 14, 727--733.

\bibitem{lss} S.~V.~Shadrin, I.~I.~Shneiberg, Belorousski-Pandharipande relation in dGBV algebras, math.QA/0507107.

\bibitem{t} U.~Tillmann, Vanishing of the Batalin-Vilkovisky algebra structure for TCFTs,
Comm. Math. Phys. 205 (1999), no. 2, 283--286.

\bibitem{w} E.~Witten, Two-dimensional gravity and intersection theory on moduli space, Surveys in Differential Geometry (Cambridge, Mass, 1990), vol. 1, Lehigh University, Pennsylvania, 1991, pp. 243--269.

\end{thebibliography}
\end{document}